\documentclass[12pt]{article}
\usepackage{amsmath,amssymb,amsthm}

\newtheorem{lemma}{Lemma}[section]

\newtheorem{theorem}{Theorem}[section]
\newtheorem{corollary}{Corollary}[section]

\theoremstyle{definition}

\begin{document}
\title{\LARGE \bf On equi-Weyl almost periodic selections of multivalued maps}
\author{\Large L.I.~Danilov \medskip \\
\large Physical-Technical Institute \\
\large Russia, 426000, Izhevsk, Kirov st., 132 \\
\large e-mail: danilov@otf.pti.udm.ru}
\date{}
\maketitle

\begin{abstract}
We prove that equi-Weyl almost periodic multivalued maps ${\bf R}\ni t\to
F(t)\in {\rm cl}\, {\mathcal U}$ have equi-Weyl almost periodic selections, where
${\rm cl}\, {\mathcal U}$ is the collection of non-empty closed sets of a 
complete metric space ${\mathcal U}$. \par
\end{abstract}

\large

\section*{Introduction}
Let $({\mathcal U},\rho )$ be a complete metric space and let $({\rm cl}_b \,
{\mathcal U},{\rm dist})$ be the metric space of non-empty closed bounded sets
$A\subset {\mathcal U}$ with the Hausdorff metric
$$
{\rm dist}(A,B)={\rm dist}_{\rho}(A,B)=\max \, \bigl\{ \sup\limits_{x\in A} \rho
(x,B),\sup\limits_{x\in B} \rho (x,A)\bigr\} \, ,\ A,B\in {\rm cl}_b \, {\mathcal 
U}\, ,
$$
where $\rho (x,F)=\inf_{y\in F}\rho (x,y)$ is a distance from a point 
$x\in {\mathcal U}$ to a non-empty set $F\subset {\mathcal U}$. The metric space
$({\rm cl}_b \, {\mathcal U},{\rm dist})$ is also complete. Let ${\rm cl}\, 
{\mathcal U}$ be the collection of non-empty closed sets $A\subset {\mathcal U}$. 
The closure of a set $A\subset {\mathcal U}$ will be denoted by $\overline A$.

In this paper we prove the existence of equi-Weyl almost periodic (a.p.) 
selections of equi-Weyl a.p. multivalued maps ${\bf R}\ni t\to F(t)\in {\rm cl}
\, {\mathcal U}$. These results are used in the study of a.p. solutions of
differential inclusions \cite{A,ABL}.

It is known that the Bohr a.p. multivalued maps ${\bf R}\ni t\to F(t)\in {\rm 
cl}_b \, {\mathcal U}$ need not have Bohr a.p. selections. In particular, in a
finite-dimensional Euclidean space it is not always possible to supplement a
Bohr a.p. orthogonal frame to form a Bohr a.p. orthogonal basis \cite{BVLL}. It
is also not hard to construct an example of continuous function ${\bf R}\ni t\to 
F(t)\in {\bf R}^2\backslash \{ 0\} $ which is not Bohr a.p. such that ${\bf R}\ni 
t\to \{ f(t),-f(t)\} \in {\rm cl}_b \, {\bf R}^2$ is a two-valued Bohr a.p. map
\cite{D1}.

In \cite{DSh}, in the case of a separable Banach space ${\mathcal H}$ it has been
proved that for any Stepanov a.p. (of degree 1) multivalued map $F\in S_1({\bf 
R},{\rm cl}_b {\mathcal H})$ there exists a countable set of Stepanov a.p. (of
degree 1) selections $f_j\in S_1({\bf R},{\mathcal H})$, $j\in {\bf N}$, such 
that ${\rm Mod}\, f_j\subset {\rm Mod}\, F$ and $F(t)=\overline {\bigcup_jf_j
(t)}$ for almost every (a.e.) $t\in {\bf R}$. The Stepanov a.p. selections of
multivalued maps ${\bf R}\ni t\to F(t)\in {\rm cl}\, {\mathcal U}$ have been studied
in \cite{D1}, \cite{DI} -- \cite{D6}. In \cite{D6}, in particular, for a compact
metric space ${\mathcal U}$ it has been proved that a multivalued map ${\bf R}\ni t\to 
F(t)\in {\rm cl}_b \, {\mathcal U}$ belongs to the space $S_1({\bf R},{\rm cl}_b 
\, {\mathcal U})$ if and only if there exists a countable set of selections $f_j
\in S_1({\bf R},{\mathcal U})$, $j\in {\bf N}$, for which $F(t)=\overline 
{\bigcup_jf_j(t)}$ almost everywhere (a.e) and $\{ f_j(.):j\in {\bf N}\} $ is a
precompact set in $L_{\infty}({\bf R},{\mathcal U})$ (furthermore, the selections 
$f_j\in S_1({\bf R},{\mathcal U})$ of the multivalued map $F(.)\in S_1({\bf R},
{\rm cl}_b \, {\mathcal U})$ can be chosen such that ${\rm Mod}\, f_j\subset {\rm 
Mod}\, F$).

In the proofs suggested in this paper we use the modifications for equi-Weyl case
of some results from \cite{D1} and \cite{D2}.

In Section 1 we give definitions and formulate our main results. Some assertions,
which will be used throughout what follows, about a.p. functions are contained in
Section 2. The simple and well known assertions (see e.g. \cite{L}) are given
without proofs. In the subsequent Sections we prove the Theorems from Section 1. 

\section{Definitions and main results}
Let ${\rm meas}$ be Lebesgue measure on ${\bf R}$ and let $B_r(x)=\{ y\in 
{\mathcal U}:\rho (x,y)\leq r\} $, $x\in {\mathcal U}$, $r\geq 0$. A function
$f:{\bf R}\to {\mathcal U}$ is said to be {\it elementary} if there exist points 
$x_j\in {\mathcal U}$ and disjoint measurable (in the Lebesgue sense) sets $T_j
\subset {\bf R}$, $j\in {\bf N}$, such that ${\rm meas}\, {\bf R}\backslash 
\bigcup_jT_j=0$ and $f(t)=x_j$ for all $t\in T_j\, $. We denote this function by 
$f(.)=\sum_jx_j\chi _{T_j}(.)$ (where $\chi _T(.)$ is the characteristic function
of a set $T\subset {\bf R}$). For arbitrary functions $f_j:{\bf R}\to 
{\mathcal U}$, $j\in {\bf N}$, we define the function $\sum_jf_j(.)\chi _{T_j}(.)
:{\bf R}\to {\mathcal U}$ coinciding with $f_j(.)$ on the set $T_j$ for $j\in 
{\bf N}$ (in the case of the metric space ${\mathcal U}$ the notations used are
formally incorrect, but no linear operations will be carried out on the functions
under consideration). A function $f:{\bf R}\to {\mathcal U}$ is {\it measurable}
if for any $\epsilon >0$ there exists an elementary function $f_{\epsilon}:{\bf 
R}\to {\mathcal U}$ such that
$$
{\rm ess}\, \sup\limits_{\hskip -0.7cm t\in {\bf R}}\, \rho (f(t),f_{\epsilon}(t))
<\epsilon \, .
$$
The class of measurable functions $f:{\bf R}\to{\mathcal U}$ will be denoted by
$M({\bf R},{\mathcal U})$ (functions coinciding for a.e. $t\in {\bf R}$ will be
identified). Let a point $x_0\in {\mathcal U}$ be fixed. We use the notation
$$
M_p({\bf R},{\mathcal U})\doteq \{ f\in M({\bf R},{\mathcal U}):\sup\limits_{\xi
\in {\bf R}}\, \int\limits_{\xi}^{\xi +1}{\rho }^{\, p}(f(t),x_0)\, dt< +\infty 
\} \, ,\ p\geq 1\, ,
$$
and for all $l>0$ we define the metrics on $M_p({\bf R},{\mathcal U})$:
$$
D^{(\rho )}_{p,\, l}(f,g)=\biggl( \, \sup\limits_{\xi \in {\bf R}}\, \frac 1l
\int\limits_{\xi}^{\xi +l}{\rho }^{\, p}(f(t),g(t))\, dt\biggr) ^{1/p},\ f,g\in
M_p({\bf R},{\mathcal U})\, .
$$
If $l_1\geq l$, then
$$
\biggl( \frac l{l_1}\biggr) ^{1/p}D^{(\rho )}_{p,\, l}(f,g)\leq D^{(\rho )}_{p,\,
l_1}(f,g)\leq \biggl( 1+\frac l{l_1}\biggr) ^{1/p}D^{(\rho )}_{p,\, l}(f,g)\, ,
$$
therefore, the metrics $D^{(\rho )}_{p,\, l}\, $, $l>0$, are equivalent and 
there exists the limit
$$
D^{(\rho ),W}_p(f,g)=\lim\limits_{l\to +\infty }\, D^{(\rho )}_{p,\, l}(f,g)=
\inf\limits_{l>0}\, D^{(\rho )}_{p,\, l}(f,g)\, ,\ f,g\in M_p({\bf R},{\mathcal 
U})\, .
$$

For a Banach space ${\mathcal U}=({\mathcal H},\| .\| )$ ($\rho (x,y)=
\| x-y\| $, $x,y\in {\mathcal H}$) we denote by
$$
\| f\| _{p,\, l}=\biggl( \, \sup\limits_{\xi \in {\bf R}}\, \frac 1l
\int\limits_{\xi}^{\xi +l}\| f(t)\| ^p dt\biggr) ^{1/p},\ l>0\, ,
$$
and by
$$
\| f\| _p=\lim\limits_{l\to +\infty }\, \| f\| _{p,\, l}
$$
the norms and the seminorm on the linear space $M_p({\bf R},{\mathcal H})\ni f$, 
$p\geq 1$. In what follows, it is convenient to assume the Banach space 
$({\mathcal H},\| .\| )$ to be complex. If the Banach space ${\mathcal H}$ is
real, then we can consider the complexification ${\mathcal H}+i{\mathcal H}$
identifying the space ${\mathcal H}$ with the real subspace (the norm $\| .\| 
_{{\mathcal H}+i{\mathcal H}}$ on the real subspace coincides with the norm
$\| .\| $).

A set $T\subset {\bf R}$ is called {\it relatively dense} if there exists a
number $a>0$ such that $[\xi ,\xi +a]\cap T\neq \emptyset $ for all $\xi \in 
{\bf R}$. A continuous function $f\in C({\bf R},{\mathcal U})$ belongs to the
space $CAP\, ({\bf R},{\mathcal U})$ of {\it Bohr a.p.} functions if for any 
$\epsilon >0$ there exists a relatively dense set $T\subset {\bf R}$ such that
the inequality
$$
\sup\limits_{t\in {\bf R}}\, \rho (f(t),f(t+\tau ))<\epsilon 
$$
holds for all $\tau \in T$. A number $\tau \in {\bf R}$ is called an {\it 
$(\epsilon ,D^{(\rho )}_{p,\, l})$-almost period} of function $f\in M_p({\bf 
R},{\mathcal U})$, $\epsilon >0$, if $D^{(\rho )}_{p,\, l}(f(.),f(.+\tau ))<
\epsilon $. A function $f\in M_p({\bf R},{\mathcal U})$, $p\geq 1$, belongs to
the space $S_p({\bf R},{\mathcal U})$ of {\it Stepanov a.p.} functions {\it of
degree $p$} if for any $\epsilon >0$ the set of $(\epsilon ,D^{(\rho )}_{p,\, 
1})$-almost periods of $f$ is relatively dense. A function $f\in M_p({\bf R},
{\mathcal U})$, $p\geq 1$, belongs to the space $W_p({\bf R},{\mathcal U})$ of
{\it equi-Weyl a.p.} functions {\it of degree $p$} if for any $\epsilon >0$ 
there exists a number $l=l(\epsilon ,f)>0$ such that the set of $(\epsilon ,
D^{(\rho )}_{p,\, l})$-almost periods of $f$ is relatively dense. We have the
inclusions  $CAP\, ({\bf R},{\mathcal U})\subset S_p({\bf R},{\mathcal U})
\subset W_p({\bf R},{\mathcal U})$.

The a.p. functions $f\in W_p({\bf R},{\mathcal U})$ are called equi-Weyl to
distinguish them from Weyl a.p. functions (a function $f\in M_p({\bf R},
{\mathcal U})$ is said to be Weyl a.p. function if for any $\epsilon >0$
there is a relatively dense set $T\subset {\bf R}$ such that $D^{(\rho ),W}
_p(f(.),f(.+\tau ))<\epsilon $ for all $\tau \in T$). In \cite{L}, the
functions $f\in W_p({\bf R},{\mathcal U})$ are called Weyl almost periodic.

A sequence $\tau _j\in {\bf R}$, $j\in {\bf N}$, is said to be {\it 
$f$-returning} for a function $f\in W_p({\bf R},{\mathcal U})$ if for any
$\epsilon >0$ there exist numbers $l=l(\epsilon ,f)>0$ and $j_0\in {\bf 
N}$ such that all numbers $\tau _j:j\geq j_0$ are $(\epsilon ,D^{(\rho 
)}_{p,\, l})$-almost periods of function $f$. 

If $f\in S_p({\bf R},{\mathcal U})\subset W_p({\bf R},{\mathcal U})$, then
a sequence $\tau _j\in {\bf R}$, $j\in {\bf N}$, is $f$-returning if and only
if $D^{(\rho )}_{p,\, 1}(f(.),f(.+\tau _j))\to 0$ as $j\to +\infty $. If $f
\in CAP\, ({\bf R},{\mathcal U})\subset W_p({\bf R},{\mathcal U})$, then a
sequence $\tau _j\in {\bf R}$, $j\in {\bf N}$, is $f$-returning if and only if
$$
\sup\limits_{t\in {\bf R}}\rho (f(t),f(t+\tau _j))\to 0
$$
as $j\to +\infty $.

For a function $f\in W_p({\bf R},{\mathcal U})$ we denote by ${\rm Mod}\, f$ 
the set of numbers $\lambda \in {\bf R}$ for which $e^{\, i\lambda \tau _j}\to 1$
($i^2=-1$) as $j\to +\infty $ for all $f$-returning sequences $\tau _j\, $. The
set ${\rm Mod}\, f$ is a module (additive group) in ${\bf R}$. If $D^{(\rho ),W}
_p(f(.),f_0(.))\neq 0$ for all constant functions $f_0(t)\equiv f_0\in {\mathcal 
U}$, $t\in {\bf R}$, then ${\rm Mod}\, f$ is a countable module (${\rm Mod}\, f 
=\{ 0\} $ otherwise).

On the space ${\mathcal U}$ we also consider the metric $\rho ^{\, \prime}(x,y)=
\min \, \{ 1,\rho (x,y)\} $, $x,y\in {\mathcal U}$. The metric space $({\mathcal 
U},\rho ^{\, \prime})$ is complete as well as $({\mathcal U},\rho )$. For all 
$f,g\in M({\bf R},{\mathcal U})=M_1({\bf R},({\mathcal U},\rho ^{\, \prime}))$ 
we use the notations
$$
D^{(\rho )}_l(f,g)=D^{(\rho ^{\, \prime})}_{1,\, l}(f,g)= 
\sup\limits_{\xi \in {\bf R}}\, \frac 1l \int\limits_{\xi}^{\xi +l}{\rho }^{\, 
\prime}(f(t),g(t))\, dt\, ,\ l>0\, ,
$$ $$
D^{(\rho ),W}(f,g)=\lim\limits_{l\to +\infty}D^{(\rho )}_l(f,g)\, .
$$
Let $S({\bf R},{\mathcal U})\doteq S_1({\bf R},({\mathcal U},\rho ^{\, 
\prime}))$, $W({\bf R},{\mathcal U})\doteq W_1({\bf R},({\mathcal U},\rho ^{\, 
\prime}))$. We have $S({\bf R},{\mathcal U})\subset W({\bf R},{\mathcal U})$, 
$W_p({\bf R},{\mathcal U})\subset W({\bf R},{\mathcal U})$, $p\geq 1$. A
sequence $\tau _j\in {\bf R}$, $j\in {\bf N}$, is $f$-returning for a function
$f\in W_p({\bf R},{\mathcal U})$ if and only if it is $f$-returning for a
function $f$ considered as the element of the space $W({\bf R},{\mathcal U})$ 
(the set of $f$-returning sequences is determined only by the a.p. function
itself and does not depend on the spaces under consideration of a.p. functions
that include the function $f$).

We shall denote by $W({\bf R})$ the collection of measurable (in the Lebesgue
sense) sets $T\subset {\bf R}$ such that $\chi _T\in W_1({\bf R},{\bf R})$. For a
set $T\in W({\bf R})$ let ${\rm Mod}\, T={\rm Mod}\, \chi _T\, $.

For an arbitrary module (additive group) $\Lambda \subset {\bf R}$ let
${\mathfrak M}^{\, (W)}(\Lambda )$ be the set of sequences $\{ T_j\} _{j\in {\bf 
N}}$ of disjoint sets $T_j\in W({\bf R})$ such that ${\rm Mod}\, T_j\subset 
\Lambda $, ${\rm meas}\, {\bf R}\backslash \bigcup_jT_j=0$ and $\| \chi _{{\bf 
R}\backslash \bigcup_{j\leq N}T_j}(.)\| _1\to 0$ as $N\to +\infty $. We shall
assume that the set ${\mathfrak M}^{\, (W)}(\Lambda )$ includes the corresponding
finite sequences $\{ T_j\} _{j=1,\dots ,N}\, $ as well, which can always be
supplemented by empty sets to form denumerable sequences. The sets $T_j$ of
sequences $\{ T_j\} \in {\mathfrak M}^{\, (W)}(\Lambda )$ will also be enumerated
by means of several indices.

If $\Lambda _j\subset {\bf R}$ are arbitrary modules, then by $\sum_j 
\Lambda _j $ (or by $\Lambda _1+\dots +\Lambda _n$ for finitely many modules 
$\Lambda _j\, $, $j=1,\dots ,n$) we denote the sum of modules, that is, the
smallest module (additive group) in ${\bf R}$ containing all the sets $\Lambda 
_j\, $.

\begin{theorem} \label{th1}
Suppose that $\{ T_j\} \in {\mathfrak M}^{\, (W)}({\bf R})$ and $f_j\in W({\bf 
R},{\mathcal U})$ for $j\in {\bf N}$. Then
$$
\sum\limits_jf_j(.)\chi _{T_j}(.)\in W({\bf R},{\mathcal U})
$$
and
$$
{\rm Mod}\, \sum\limits_jf_j(.)\chi _{T_j}(.)\subset \sum\limits_j{\rm Mod}\, f_j
+\sum\limits_j{\rm Mod}\, T_j\, .  \eqno (1)
$$
\end{theorem}

{\bf Remark 1}. Under the assumptions of Theorem \ref{th1} for indices $j\in {\bf 
N}$ such that $\| \chi _{T_j}\| _1=0$ (in this case ${\rm Mod}\, T_j=\{ 0\} $) we
can choose arbitrary functions $f_j\in M({\bf R},{\mathcal U})$ and delete
these indices in the summation on the right-hand side of inclusion (1).
\vskip 0.2cm

The following Theorem is analogous to the Theorem on uniform approximation of
Stepanov a.p. functions \cite{D1,D4,D5} and plays a key role in this paper.

\begin{theorem} \label{th2}
Let $f\in W({\bf R},{\mathcal U})$. Then for any $\epsilon >0$ there exist a
sequence $\{ T_j\} \in {\mathfrak M}^{\, (W)}({\rm Mod}\, f)$ and points $x_j\in 
{\mathcal U}$, $j\in {\bf N}$, such that $\rho (f(t),x_j)<\epsilon $ for all $t
\in T_j\, $, $j\in {\bf N}$.
\end{theorem}

Theorem \ref{th2} is proved in Section 3.

Let ${\rm dist}_{\rho ^{\, \prime}}$ be the Hausdorff metric on ${\rm cl}\, 
{\mathcal U}={\rm cl}_b \, ({\mathcal U},\rho ^{\, \prime})$ corresponding to the 
metric 
$\rho ^{\, \prime}$. The metric space $({\rm cl}\, {\mathcal U},{\rm dist}_{\rho 
^{\, \prime}})$ is complete. Since  ${\rm dist}^{\, \prime}(A,B)\doteq \min \, \{ 
1,{\rm dist}\, (A,B)\}={\rm dist}_{\rho ^{\, \prime}}(A,B)$ for all $A,B\in {\rm 
cl}_b \, {\mathcal U}$, it follows that the embedding $({\rm cl}_b\, {\mathcal 
U},{\rm dist}^{\, \prime})\subset ({\rm cl}\, {\mathcal U},{\rm dist}_{\rho ^{\, 
\prime}})$ is isometric. We define the spaces $W({\bf R},{\rm cl}_b \, {\mathcal 
U})$ and $W_p({\bf R},{\rm cl}_b \, {\mathcal U})$, $p\geq 1$, of {\it equi-Weyl
a.p. multivalued maps} ${\bf R}\ni t\to F(t)\in {\rm cl}_b\, {\mathcal U}$ as the 
spaces of equi-Weyl a.p. functions taking values in the metric space $({\rm cl}_b 
\, {\mathcal U},{\rm dist})$. Let $W({\bf R},{\rm cl}\, {\mathcal U})\doteq W_1
({\bf R},({\rm cl}\, {\mathcal U},{\rm dist}_{\rho ^{\, \prime}}))$. The 
following inclusions $W_p({\bf R},{\rm cl}_b \, {\mathcal U})\subset W_1({\bf R},
{\rm cl}_b\, {\mathcal U})\subset W({\bf R},{\rm cl}_b\, {\mathcal U})\subset 
W({\bf R},{\rm cl}\, {\mathcal U})$ hold.

Let us denote by ${\mathcal N}$ the set of non-decreasing functions $[0,+\infty )
\ni t\to \eta (t)\in {\bf R}$ such that $\eta (0)=0$ and $\eta (t)>0$ for all 
$t>0$.

\begin{theorem} \label{th3}
Let $({\mathcal U},\rho )$ be a complete metric space, let $F\in W({\bf R},{\rm 
cl}\, {\mathcal U})$ and let $g\in W({\bf R},{\mathcal U})$. Then for any
function $\eta \in {\mathcal N}$ there exists a function $f\in W({\bf R},
{\mathcal U})$ such that ${\rm Mod}\, f\subset {\rm Mod}\, F+{\rm Mod}\, g$, 
$f(t)\in F(t)$ a.e. and $\rho (f(t),g(t))\leq \rho (g(t),F(t))+\eta (\rho (g(t),
F(t)))$ a.e. Moreover, if $F\in W_p({\bf R},{\rm cl}_b\, {\mathcal U})\subset 
W({\bf R},{\rm cl}\, {\mathcal U})$, $p\geq 1$, then $f\in W_p({\bf R},{\mathcal 
U})$.
\end{theorem}

Theorem \ref{th3} is the main result of the paper on equi-Weyl a.p. selections
of multivalued maps and is proved in Section 5.

\section{Some properties of equi-Weyl a.p. functions}

For functions $f,g\in M_p({\bf R},{\mathcal U})$, $p\geq 1$, we set
$$
J_p(f,g)=\lim\limits_{\delta \to +0}\, \lim\limits_{l_0\to +\infty }\,
\sup\limits_{l\geq l_0}\, \sup\limits_{\xi \in {\bf R}}\, \biggl( \, \frac 1l\,
\sup\limits_{T\subset [\xi ,\xi +l]\, :\, {\rm meas}\, T\leq \delta l}\ 
\int\limits_T\rho ^{\, p}(f(t),g(t))\, dt\biggr) ^{1/p};
$$ $$
J_p(f,g)\leq \lim\limits_{l_0\to +\infty }\, \sup\limits_{l\geq l_0}\, D^{(\rho
)}_{p,\, l}(f,g)=D^{(\rho ),W}_p(f,g)\, .
$$
We use the notation 
$$
M_p^{\sharp}({\bf R},{\mathcal U})=\{ f\in M_p({\bf R},{\mathcal U}): J_p(f(.),
x_0(.))=0\} \, ,
$$
where $x_0(t)\equiv x_0\, $, $t\in {\bf R}$. The set $M_p^{\sharp}({\bf R},
{\mathcal U})$ does not depend on the choice of the point $x_0\in {\mathcal U}$.

The following simple Lemmas are used in the proof of Theorem \ref{th4}.

\begin{lemma} \label{l1}
A function $f\in M({\bf R},{\mathcal U})$ belongs to the space $W({\bf R},
{\mathcal U})$ of equi-Weyl a.p. functions if and only if for any $\epsilon ,
\delta >0$ there exists a number $l>0$ such that the inequality 
$$
\sup\limits_{\xi \in {\bf R}}\, {\rm meas}\, \{ t\in [\xi ,\xi +l]: \rho (f(t),
f(t+\tau ))\geq \delta \} <\epsilon l\, .
$$
holds for relatively dense set of numbers $\tau \in {\bf R}$. Furthermore, the
sequence $\tau _j\in {\bf R}$, $j\in {\bf N}$, is $f$-returning for a function 
$f\in W({\bf R},{\mathcal U})$ if and only if for any $\epsilon ,\delta >0$ there
exist numbers $l>0$ and $j_0\in {\bf N}$ such that for all $j\in {\bf N}:j\geq 
j_0$ 
$$
\sup\limits_{\xi \in {\bf R}}\, {\rm meas}\, \{ t\in [\xi ,\xi +l]: \rho (f(t),
f(t+\tau _j))\geq \delta \} <\epsilon l\, .
$$
\end{lemma}

For a function $f\in M({\bf R},{\mathcal U})$ and a number $R\geq 0$ we define
the function
$$
{\bf R}\ni t\to f_R(x_0;t)=\left\{
\begin{array}{ll}
f(t)  &{\text {if}}\ \rho (f(t),x_0)\leq R\, ,\\ [0.2cm]
x_0   &{\text {if}}\ \rho (f(t),x_0)>R\, .
\end{array}
\right.
$$

\begin{lemma} \label{l2}
Let $f\in W_p({\bf R},{\mathcal U})$, $p\geq 1$, $x_0\in {\mathcal U}$. Then 
$D^{(\rho ),W}_p(f(.),f_R(x_0;.))\to 0$ as $R\to +\infty $.
\end{lemma}

\begin{theorem} \label{th4}
For all $p\geq 1$ 
$$
W_p({\bf R},{\mathcal U})=W({\bf R},{\mathcal U})\bigcap M_p^{\sharp}({\bf R},
{\mathcal U})\, .
$$
\end{theorem}

\begin{proof} 
Let $f\in W_p({\bf R},{\mathcal U})$. By Lemma \ref{l2} (in view of
boundedness of functions $f_R(x_0;.)$), 
$$
J_p(f(.),x_0(.))\leq J_p(f(.),f_R(x_0;.))+J_p(f_R(x_0;.),x_0))=
$$ $$
=J_p(f(.),f_R(x_0;.))\leq D^{(\rho ),W}_p(f(.),f_R(x_0;.))\to 0
$$
as $R\to +\infty $. Hence $W_p({\bf R},{\mathcal U})\subset M_p^{\sharp}({\bf 
R},{\mathcal U})$, and therefore, $W_p({\bf R},{\mathcal U})\subset W({\bf R},
{\mathcal U})\bigcap M_p^{\sharp}({\bf R},{\mathcal U})$. To prove the reverse
inclusion it is necessary to use Lemma \ref{l1}.
\end{proof}

Let $({\mathcal H},\| .\|)$ be a complex Banach space. For any function $f\in 
W_p({\bf R},{\mathcal H})$ and any $\lambda \in {\bf R}$ there exists a limit
$$
\lim\limits_{a\to +\infty }\, \frac 1{2a}\, \int\limits_{-a}^ae^{-i\lambda t}f(t)
\, dt=M\, \{ e^{-i\lambda t}f\}
$$
(the integral is defined in the sense of Bochner). We denote by $\Lambda \{ f\} $
the set of Fourier exponents of a function $f\in W_p({\bf R},{\mathcal H}):
\Lambda \{ f\} =\{ \lambda \in {\bf R}:M\, \{ e^{-i\lambda t}f\} \neq 0\} $. Let
${\rm Mod}\, (\Lambda \{ f\} )$ be {\it the module of the Fourier exponents} of a
function $f$ (that is, the smallest additive group in ${\bf R}$ containing the set
$\Lambda \{ f\} $).

\begin{lemma} \label{l3}
Let $f\in W_p({\bf R},{\mathcal U})$. Then for any sequence $\tau _j\in {\bf R}$, 
$j\in {\bf N}$, the following three conditions are equivalent:

(1) $\{ \tau _j\} $ is $f$-returning sequence;

(2) $D^{(\rho ),W}_p(f(.),f(.+\tau _j))\to 0$ as $j\to +\infty $,

(3) $e^{\, i\lambda \tau _j}\to 1$ as $j\to +\infty $ for all $\lambda \in
{\rm Mod}\, f$.\\
If ${\mathcal U}=({\mathcal H},\| .\|)$, then ${\rm Mod}\, f={\rm Mod}\, (\Lambda
\{ f\} )$. 
\end{lemma}

Suppose that $f\in W({\bf R},{\mathcal U})$ and $f_j\in W({\bf R},{\mathcal 
U}_j)$, $j\in {\bf N}$, where the ${\mathcal U}_j$ are (complete) metric spaces.
Then ${\rm Mod}\, f\subset \sum_j{\rm Mod}\, f_j$ if and only if every sequence
$\tau _k\in {\bf R}$, $k\in {\bf N}$, that is $f_j$-returning for all $j\in {\bf 
N}$ is $f$-returning. In particular, for the case when $f_j\in W({\bf R},
{\mathcal U}_j)$, $j=1,2$, we have ${\rm Mod}\, f_1\subset {\rm Mod}\, f_2$ if 
and only if every $f_2$-returning sequence $\{\tau _k\}$ is $f_1$-returning.

\begin{lemma} \label{l4}
Suppose that $f\in M_p({\bf R},{\mathcal U})$, $f_j\in W_p({\bf R},{\mathcal 
U})$, $j\in {\bf N}$, and \\ $D^{(\rho ),W}_p(f,f_j)\to 0$ as $j\to +\infty $. 
Then $f\in W_p({\bf R},{\mathcal U})$ and ${\rm Mod}\, f\subset \sum_j{\rm Mod}\, 
f_j\, $.
\end{lemma}

\begin{corollary} \label{c1}
Suppose that $f\in M({\bf R},{\mathcal U})$, $f_j\in W({\bf R},{\mathcal U})$, 
$j\in {\bf N}$, and $D^{(\rho ),W}(f,f_j)\to 0$ as $j\to +\infty $. Then $f\in 
W({\bf R},{\mathcal U})$ and ${\rm Mod}\, f\subset \sum_j{\rm Mod}\, f_j\, $.
\end{corollary}

\begin{theorem}[see e.g. \cite{L}] \label{th5}
Let $({\mathcal H},\| .\| )$ be a complex Banach space and let $f\in W_p({\bf R},
{\mathcal H})$, $p\geq 1$. Then for any $\epsilon >0$ there is a trigonometric
polynomial 
$$
f_{\epsilon}(t)=\sum\limits_{j=1}^{N(\epsilon )}c^{(\epsilon )}_je^{\, i{\lambda}
^{(\epsilon )}_jt}\, ,\ t\in {\bf R}\, ,
$$ 
where $c^{(\epsilon )}_j\in {\mathcal H}$, ${\lambda}^{(\epsilon )}_j\in {\bf R}$
(and the sum contains only a finite number of terms), such that $\| f-f
_{\epsilon}\| _p<\epsilon $ and $\Lambda \{ f_{\epsilon }\} \subset \Lambda \{ 
f\} $.
\end{theorem}

The following Theorem is a consequence of Theorem \ref{th5}.

\begin{theorem} \label{th6}
Let $f\in W({\bf R},{\mathcal H})$. Then for any $\epsilon >0$ there is a
trigonometric polynomial $f_{\epsilon}\in CAP\, ({\bf R},{\mathcal H})$ such that 
$D^{(\rho ),W}(f,f_{\epsilon })<\epsilon $ and ${\rm Mod}\, f_{\epsilon}\subset 
{\rm Mod}\, f$.
\end{theorem}

\begin{corollary} \label{c2}
Let $f_1\, ,f_2\in W({\bf R},{\mathcal H})$, then $f_1+f_2\in W({\bf R},{\mathcal 
H})$ and ${\rm Mod}\, (f_1+f_2)\subset {\rm Mod}\, f_1+{\rm Mod}\, f_2\, $. If 
$f\in W({\bf R},{\mathcal H})$, $g\in W({\bf R},{\bf C})$, then also $gf\in 
W({\bf R},{\mathcal H})$ and ${\rm Mod}\, gf\subset {\rm Mod}\, f+{\rm Mod}\, g$.
\end{corollary}

For a set $T\in W({\bf R})$ we have ${\bf R}\backslash T\in W({\bf R})$ and ${\rm
Mod}\, {\bf R}\backslash T={\rm Mod}\, T$.

\begin{lemma} \label{l5}
Let $T_1\, ,T_2\in W({\bf R})$. Then $T_1\bigcup T_2\in W({\bf R})$, $T_1\bigcap 
T_2\in W({\bf R})$, $T_1\backslash T_2\in W({\bf R})$ and modules ${\rm Mod}\, 
T_1\bigcup T_2\, $, ${\rm Mod}\, T_1\bigcap T_2$, ${\rm Mod}\, T_1\backslash T_2 
$ are subsets (subgroups) of ${\rm Mod}\, T_1+{\rm Mod}\, T_2\, $.
\end{lemma}

\begin{corollary} \label{c3}
Let  $\Lambda $ be a module in ${\bf R}$ and let $\{ T^{(s)}_j\}_{j\in {\bf N}}
\in {\mathfrak M}^{(W)}(\Lambda )$, $s=1,2$, then also $\{ T^{(1)}_j\bigcap 
T^{(2)}_k\} _{j,k\in {\bf N}}\in {\mathfrak M}^{(W)}(\Lambda )$.
\end{corollary}

{\it Proof of Theorem \ref{th1}}. By the Fr\'echet Theorem (on isometric 
embedding of a metric space into some Banach space) \cite{LS}, we can suppose
that ${\mathcal U}=({\mathcal H},\| .\| )$. From Corollary \ref{c2} it follows
that for all $N\in {\bf N}$
$$
\sum\limits_{j=1}^Nf_j(.)\chi _{T_j}(.)\in W({\bf R},{\mathcal H})
$$
and
$$
{\rm Mod}\, \sum\limits_{j=1}^Nf_j(.)\chi _{T_j}(.)\subset \sum\limits_{j=1}^N
{\rm Mod}\, f_j+\sum\limits_{j=1}^N{\rm Mod}\, T_j\, .
$$
On the other hand, we have
$$
D^{(\rho ),W}\bigl( \, \sum\limits_{j=1}^{+\infty }f_j(.)\chi _{T_j}(.)\, ,\, 
\sum\limits_{j=1}^Nf_j(.)\chi _{T_j}(.)\bigr) \to 0
$$
as $N\to +\infty $. To complete the proof it remains to apply Corollary \ref{c1}.
\hfill  $\square$

\section{Proof of Theorem \ref{th2}}

For $h\in ({\mathcal H},\| .\| )$ we set
$$
{\rm sgn}\, h=\left\{
\begin{array}{ll}
\frac h{\| h\| }\ \  &{\text {if}}\ h\neq 0\, ,\\ [0.2cm]
0\ \  &{\text {if}}\ h=0\, .
\end{array}
\right.
$$

\begin{lemma} \label{l6}
Let $f\in W({\bf R},{\mathcal H})$. Suppose that for any $\epsilon >0$ there are
numbers $\delta >0$ and $l>0$ such that
$$
\sup\limits_{\xi \in {\bf R}}\, {\rm meas}\, \{ t\in [\xi ,\xi +l]:\| f(t)\| <
\delta \} <\epsilon l\, .  \eqno (2)
$$
Then ${\rm sgn}\, f(.)\in W_1({\bf R},{\mathcal H})$ and ${\rm Mod}\, {\rm sgn}\,
f(.)\subset {\rm Mod}\, f(.)$ (furthermore, for the set $T=\{ t\in {\bf R}:f(t)=0
\} $ we have $\| \chi _T(.)\| _1=0$).
\end{lemma}

\begin{proof} For all $j\in {\bf N}$ define the functions
$$
{\mathcal H}\ni h\to {\mathcal F}_j(h)=\left\{
\begin{array}{ll}
jh\ \  &{\text {if}}\ \| h\| \leq \frac 1j\, ,\\ [0.2cm]
\frac h{\| h\| }\ \  &{\text {if}}\ \| h\| >\frac 1j\, .
\end{array}
\right.
$$
These functions are uniformly continuous and bounded, therefore 
${\mathcal F}_j(f(.))\in W_1({\bf R},{\mathcal H})$ and ${\rm Mod}\, {\mathcal 
F}_j(f(.))\subset {\rm Mod}\, f(.)$. From the condition (2) it follows that 
$\| \chi _T(.)\| _1=0$ and $\| {\rm sgn}\, f(.)-{\mathcal F}_j(f(.))\| _1\to 0$ 
as $j\to +\infty $. Hence (in view of Lemma \ref{l4}) ${\rm sgn}\, f(.)\in W_1
({\bf R},{\mathcal H})$ and ${\rm Mod}\, {\rm sgn}\, f(.)\subset \sum_j
{\rm Mod}\, {\mathcal F}_j(f(.))\subset {\rm Mod}\, f(.)$. 
\end{proof}

\begin{lemma} \label{l7} 
Let $f\in W({\bf R},{\mathcal U})$. Then for any $\epsilon ,\delta >0$ there 
exist a number $l>0$ and finitely many points $x_j\in {\mathcal U}$, $j=1,\dots 
,N$, such that
$$
\sup\limits_{\xi \in {\bf R}}\, {\rm meas}\, \{ t\in [\xi ,\xi +l]:\rho (f(t),\,
\bigcup\limits_{j=1}^Nx_j)\geq \delta \} <\epsilon l\, .
$$
\end{lemma}

\begin{proof} 
Lemma \ref{l1} implies that for any $\epsilon ,\delta >0$ there exist a number 
$l>0$ and a relatively dense set $T\subset {\bf R}$ such that the inequality
$$
\sup\limits_{\xi \in {\bf R}}\, {\rm meas}\, \{ t\in [\xi ,\xi +l]:\rho (f(t),
f(t+\tau ))\geq \frac {\delta}2\, \} <\frac {\epsilon}2\ l\, 
$$
holds for all $\tau \in T$. We choose a number $a\geq \frac l2$ such that for any 
$\xi \in {\bf R}$ there is a number $\tau (\xi )\in T$ for which $[\xi +\tau (\xi 
),\xi +l+\tau (\xi )]\subset [-a,a]$. Since the function $f$ is measurable, we 
have 
$$
{\rm meas}\, \{ t\in [-a,a]:\rho (f(t),\, \bigcup\limits_{j=1}^Nx_j)\geq \frac
{\delta}2\, \} <\frac {\epsilon}2\ l\, 
$$
for some finite set of points $x_j\in {\mathcal U}$, $j=1,\dots ,N$, and hence
$$
\sup\limits_{\xi \in {\bf R}}\, {\rm meas}\, \{ t\in [\xi ,\xi +l]:\rho (f(t),\,
\bigcup\limits_{j=1}^Nx_j)\geq \delta \} \leq
$$ $$
\leq \sup\limits_{\xi \in {\bf R}}\ \biggl( {\rm meas}\, \{ t\in [\xi ,\xi +l]:
\rho (f(t),f(t+\tau (\xi )))\geq \frac {\delta}2\, \} \, +
$$ $$
+\, {\rm meas}\, \{ t\in [\xi ,\xi +l]:\rho (f(t+\tau (\xi )),\bigcup\limits_{j=1}
^Nx_j)\geq \frac {\delta}2\, \} \biggr) \leq
$$ $$
\leq \frac {\epsilon}2\ l+{\rm meas}\, \{ t\in [-a,a]:\rho (f(t),\, 
\bigcup\limits_{j=1}^Nx_j)\geq \frac {\delta}2\, \} <\frac {\epsilon}2\ l+\frac 
{\epsilon}2\ l=\epsilon l\, . 
$$
\end {proof}

\begin{corollary} \label{c4}
Let $f\in W({\bf R},{\mathcal U})$. Then for any $\delta >0$ there are points 
$x_j\in {\mathcal U}$, $j\in {\bf N}$, such that

(1) for all $a>0$
$$
\lim\limits_{N\to +\infty }\, {\rm meas}\, \{ t\in [-a,a]:\rho (f(t),\, 
\bigcup\limits_{j\leq N}x_j)\geq \delta \} =0\, ,
$$

(2) for any $\epsilon >0$ there exist numbers $l>0$ and $N\in {\bf N}$ such that
$$
\sup\limits_{\xi \in {\bf R}}\, {\rm meas}\, \{ t\in [\xi ,\xi +l]:\rho (f(t),\, 
\bigcup\limits_{j\leq N}x_j)\geq \delta \} <\epsilon l\, .
$$
\end{corollary}

Let ${\mathcal A}^{(W)}$ be the collection of sets $\mathbb F\subset W({\bf R},
{\bf R})$ such that for any $\epsilon >0$ there exist numbers $l=l(\epsilon ,
{\mathbb F})>0$ and $\tau _0=\tau _0(\epsilon ,{\mathbb F})>0$ for which
$$
\sup\limits_{f\in {\mathbb F}}\, \sup\limits_{\tau \in [0,\tau _0]}\, D^{(\rho )}
_l(f(.),f(.+\tau ))<\epsilon
$$
($\rho (x,y)=|x-y|$, $x,y\in {\bf R}$).

\begin{lemma}[\cite{L}] \label{l8}
Let $f\in W_p({\bf R},{\mathcal U})$. Then for any $\epsilon >0$ there are 
numbers $l>0$ and $\tau _0>0$ such that the inequality $D^{(\rho )}_{p,\, l}(f(.)
,f(.+\tau ))<\epsilon $ holds for all $\tau \in [0,\tau _0]$.
\end{lemma}

From Lemma \ref{l8} it follows that $\{ f\} \in {\mathcal A}^{(W)}$ for any 
function $f\in W({\bf R},{\bf R})$.

The following Theorem is proved in Section 4 and its special case for the set 
${\mathbb F}=\{ f\} $, $f\in W({\bf R},{\bf R})$, is essentially used in the 
proof of Theorem \ref{th2}.

\begin{theorem} \label{th7}
Let ${\mathbb F}\in {\mathcal A}^{(W)}$, $\Delta >0$, $b>0$, $\epsilon \in (0,1]
$. Then there exist $b$-periodic function $g(.)\in C({\bf R},{\bf R})$ dependent
on ${\mathbb F}$, $\Delta $, $b$, but not on the number $\epsilon $, for which 
$\| g\| _{L_{\infty}({\bf R},{\bf R})}<\Delta $, and numbers $\delta = \delta 
(\epsilon ,\Delta )>0$, $l=l(\epsilon ,\Delta ,{\mathbb F})>0$ such that for all
functions $f\in {\mathbb F}$
$$
\sup\limits_{\xi \in {\bf R}}\, {\rm meas}\, \{ t\in [\xi ,\xi +l]:|f(t)+g(t)|<
\delta \} <\epsilon l\, .
$$
\end{theorem}

\begin{corollary} \label{c5}
Let $f\in W({\bf R},{\bf R})$. Then for any $a\in {\bf R}$ and $\epsilon >0$ 
there is a set $T\in W({\bf R})$ such that ${\rm Mod}\, T\subset {\rm Mod}\, f$, 
$f(t)<a+\epsilon $ for all $t\in T$ and $f(t)>a$ for a.e. $t\in {\bf R}\backslash 
T$. 
\end{corollary}

{\it Proof of Theorem \ref{th2}}. If ${\rm Mod}\, f=\{ 0\} $, then for some 
constant function $f_0(t)\equiv f_0\in {\mathcal U}$, $t\in {\bf R}$, we have 
$D^{(\rho ),W}(f(.),f_0(.))=0$, and therefore, there is a set $T\in W({\bf R})$ 
such that $\| \chi _T(.)\| _1=0$ and $\rho (f(t),f_0)<\epsilon $ for all $t\in 
{\bf R}\backslash T$. Next, suppose that ${\rm Mod}\, f\neq \{ 0\} $. Let
$x_j\in {\mathcal U}$, $j\in {\bf N}$, be the points defined in Corollary 
\ref{c4} for the function $f\in W({\bf R},{\mathcal U})$ and the number $\delta 
=\frac {\epsilon}3\, $. For all $j\in {\bf N}$ we have $\rho (f(.),x_j)\in W({\bf 
R},{\bf R})$ and ${\rm Mod}\, \rho (f(.),x_j)\subset {\rm Mod}\, f(.)$. We choose
a number $b>0$ such that $\frac {2\pi }b\in {\rm Mod}\, f$. Theorem \ref{th7}
implies the existence of $b$-periodic function $g_j(.)\in C({\bf R},{\bf R})$, 
$j\in {\bf N}$, such that $\| g_j\| _{L_{\infty}({\bf R},{\bf R})}<\frac 
{\epsilon }3\, $ and for any $\epsilon ^{\, \prime}\in (0,1]$ there exist numbers
$\delta _j=\delta _j(\epsilon ^{\, \prime},\epsilon )>0$ and $l_j=l_j(\epsilon 
^{\, \prime},\epsilon ,f)>0$ for which
$$
\sup\limits_{\xi \in {\bf R}}\, {\rm meas}\, \{ t\in [\xi ,\xi +l_j]:|\rho (f(t),
x_j)-\frac {2\epsilon }3+g_j(t)|<\delta _j\} < \epsilon ^{\, \prime}l_j\, .
$$
Let $T^{\, \prime}_j=\{ \, t\in {\bf R}:\rho (f(t),x_j)+g_j(t)\leq \frac {2
\epsilon }3\, \} $, $j\in {\bf N}$. According to Lemma \ref{l6}, we get $T^{\, 
\prime}_j\in W({\bf R})$ and ${\rm Mod}\, T^{\, \prime}_j\subset {\rm Mod}\, \rho 
(f(.),x_j)+\frac {2\pi }b\, {\bf Z}\subset {\rm Mod}\, f(.)$. If $t\in T^{\, \prime}_j\, $, 
then $\rho (f(t),x_j)<\epsilon $. We denote $T_1=T^{\, \prime}_1$ and $T_j=T^{\, 
\prime}_j\backslash \bigcup_{k<j}T^{\, \prime}_k$ for $j\geq 2$. The sets $T_j$,
$j\in {\bf N}$, are disjoint and $\bigcup_{j\leq N}T_j=\bigcup_{j\leq N}T^{\,
\prime}_j$ for all $N\in {\bf N}$. It follows from Lemma \ref{l5} that $T_j\in 
W({\bf R})$, ${\rm Mod}\, T_j\subset {\rm Mod}\, f$. Furthermore, $\rho (f(t),x_j)
<\epsilon $ for all $t\in T_j\, $, $j\in {\bf N}$, and for every $N\in {\bf N}$ 
and a.e. $t\in {\bf R}\backslash \bigcup_{j\leq N}T_j$ we have $\rho (f(t),x_j)
\geq \frac {\epsilon}3$ for all $j=1,\dots ,N$. Hence (see Corollary \ref{c4}) 
${\rm meas}\, {\bf R}\, \backslash \bigcup_{j}T_j=0$ and $\| \chi _{{\bf R}
\, \backslash \bigcup_jT_j}(.)\| _1\to 0$ as $N\to +\infty $, that is, 
$\{ T_j\} \in {\mathfrak M}^{(W)}({\rm Mod}\, f)$. \hfill $\square$

Theorems \ref{th1} and \ref{th2} enable us to get the following characterization 
of functions $W({\bf R},{\mathcal U})$: a function $f:{\bf R}\to {\mathcal U}$ 
belongs to $W({\bf R},{\mathcal U})$ if and only if for any $\epsilon >0$
there exist a sequence $\{ T_j\} \in {\mathfrak M}^{(W)}({\bf R})$ and points
$x_j\in {\mathcal U}$, $j\in {\bf N}$, such that the inequality $\rho (f(t),x_j)
<\epsilon $ holds for all $t\in T_j\, $, $j\in {\bf N}$.

In the following two Lemmas we consider the superposition of equi-Weyl a.p.
functions.

\begin{lemma} \label{l9}
Let $({\mathcal U},\rho )$ and $({\mathcal V},\rho _{\mathcal V})$ be complete
metric spaces, let ${\mathcal F}\in C({\mathcal U},{\mathcal V})$ and let $f\in 
W({\bf R},{\mathcal U})$. Then ${\mathcal F}(f(.))\in W({\bf R},{\mathcal V})$ 
and ${\rm Mod}\, {\mathcal F}(f(.))\subset {\rm Mod}\, f(.)$.
\end{lemma}

\begin{proof}
Let $\epsilon \in (0,1]$, $\delta >0$. By Theorem \ref{th2}, for every $k\in {\bf 
N}$ there are sequences $\{ T^{(k)}_j\} \in {\mathfrak M}^{(W)}({\rm Mod}\, f)$ 
and points $x^{(k)}_j\in {\mathcal U}$, $j\in {\bf N}$, such that $\rho (f(t),
x^{(k)}_j)<\frac 1k$ for all $t\in T^{(k)}_j$, $j\in {\bf N}$. The inequalities
$$
\| \chi _{{\bf R}\, \backslash \bigcup\limits_{j=1}^{j(k)}T^{(k)}_j}(.)\| _1<
2^{-k}\epsilon \, .
$$
hold for some numbers $j(k)\in {\bf N}$, $k\in {\bf N}$. For every $k\in {\bf N}$ 
we denote by $X_k$ the set of points $x^{(k)}_j$, $j=1,\dots ,j(k)$, for which 
for any $k^{\, \prime}=1,\dots ,\, k-1$ there exists a point $x^{(k^{\, \prime})}
_{j^{\, \prime}}$, $j^{\, \prime}=1,\dots ,j(k^{\, \prime})$, such that $\rho 
(x^{(k)}_j,x^{(k^{\, \prime})}_{j^{\, \prime}})<\frac 1k+\frac 1{k^{\, \prime}}<
\frac 2{k^{\, \prime}}\, $; $X_1=\bigcup_{j\leq j(1)}x^{(1)}_j$. Since the set 
$\bigcup_{k\in {\bf N}}X_k\subset {\mathcal U}$ is precompact and the function 
${\mathcal F}$ is continuous, it follows that there is a number $k_0\in {\bf N}$
such that for all $j_k=1,\dots ,j(k)$, where $k=1,\dots ,k_0\, $, the inequality
$\rho _{\mathcal V}({\mathcal F}(f(t)),{\mathcal F}(f(t^{\, \prime})))<\delta $
holds for all $t,t^{\, \prime}\in T^{(1)}_{j_1}\bigcap \dots \bigcap T^{(k_0)}_{j
_{k_0}}$. If $T^{(1)}_{j_1}\bigcap \dots \bigcap T^{(k_0)}_{j_{k_0}}\neq 
\emptyset $, where $j_k=1,\dots ,j(k)$, $k=1,\dots ,k_0\, $, we choose some 
numbers $t_{j_1\dots j_{k_0}}\in T^{(1)}_{j_1}\bigcap \dots \bigcap T^{(k_0)}_{j
_{k_0}}$. Let
$$
T(k_0)=\bigcap\limits_{k=1,\dots ,k_0}\, \bigcup\limits_{j_k=1}^{j(k)}T^{(k)}
_{j_k}\, .
$$
By Lemma \ref{l5} and Theorem \ref{th1},
$$
{\mathcal G}_{k_0}(.)\doteq \sum\limits_{j_k=1,\dots ,j(k);\, k=1,\dots ,k_0}
{\mathcal F}(f(t_{j_1\dots j_{k_0}}))\chi _{T^{(1)}_{j_1}\bigcap \dots \bigcap 
T^{(k_0)}_{j_{k_0}}}(.)
+y_0\chi _{{\bf R}\, \backslash T(k_0)}(.)\in W({\bf R},{\mathcal V})\, ,
$$ 
where $y_0\in {\mathcal V}$, and
$$
{\rm Mod}\, {\mathcal G}_{k_0}(.)\subset \sum\limits_{j_k=1,\dots ,j(k);\, k=1,
\dots ,k_0}{\rm Mod}\, T^{(k)}_{j_k}\subset {\rm Mod}\, f(.)\, .
$$
Furthermore, $\rho _{\mathcal V}({\mathcal F}(f(t)),{\mathcal G}_{k_0}(t))<
\delta $ for all $t\in T(k_0)$ and
$$
\| \chi _{{\bf R}\, \backslash T(k_0)}(.)\| _1\leq \sum\limits_{k=1,\dots ,k_0}\,
2^{-k}\epsilon <\epsilon \, .
$$
Hence $D^{(\rho _{\mathcal V}),W}({\mathcal F}(f(.)),{\mathcal G}_{k_0}(.))<
\epsilon +\delta $. Since the numbers $\epsilon >0$ and $\delta >0$ can be chosen
arbitraryly small, it follows from Corollary \ref{c1} that ${\mathcal F}(f(.))\in
W({\bf R},{\mathcal V})$ and ${\rm Mod}\, {\mathcal F}(f(.))\subset {\rm Mod}\,
f(.)$. 
\end{proof}

On the set $C({\mathcal U},{\mathcal V})$, where $({\mathcal U},\rho )$ and
$({\mathcal V},\rho _{\mathcal V})$ are metric spaces, we introduce the metric
$$
d_{C({\mathcal U},{\mathcal V})}({\mathcal F}_1\, ,{\mathcal F}_2)=\sup\limits
_{x\in {\mathcal U}}\, \min \, \{ 1,\rho _{\mathcal V}({\mathcal F}_1(x),
{\mathcal F}_2(x))\} \, ,\ {\mathcal F}_1\, ,\, {\mathcal F}_2\in C({\mathcal U},
{\mathcal V})\, .
$$

\begin{lemma} \label{l10} 
Let $({\mathcal U},\rho )$ and $({\mathcal V},\rho _{\mathcal V})$ be complete
metric spaces. Suppose that a function ${\bf R}\ni t\to {\mathcal F}(.;t)\in 
C({\mathcal U},{\mathcal V})$ belongs to the space $W_1({\bf R},(C({\mathcal U},
{\mathcal V}),d_{C({\mathcal U},{\mathcal V})}))$ and $f\in W({\bf R},{\mathcal 
U})$. Then ${\mathcal F}(f(.);.)\in W({\bf R},{\mathcal V})$ and ${\rm Mod}\, 
{\mathcal F}(f(.);.)\subset {\rm Mod}\, {\mathcal F}(.;.)+{\rm Mod}\, f(.)$.
\end{lemma}

\begin{proof}
Theorem \ref{th2} implies that for any $\epsilon >0$ there are a sequence $\{ 
T_j\} \in {\mathfrak M}^{(W)}({\rm Mod}\, {\mathcal F}(.;.))$ and functions 
${\mathcal F}_j\in C({\mathcal U},{\mathcal V})$, $j\in {\bf N}$, such that \\
$d_{C({\mathcal U},{\mathcal V})}({\mathcal F}(.;t),{\mathcal F}_j(.))<\epsilon $ 
for all $t\in T_j\, $, $j\in {\bf N}$. By Theorem \ref{th1} and Lemma \ref{l9},
$$
\sum\limits_{j\in {\bf N}}\, {\mathcal F}_j(f(.))\chi _{T_j}(.)\in W({\bf R},
{\mathcal V})\, ,
$$ $$
{\rm Mod}\, \sum\limits_{j\in {\bf N}}\, {\mathcal F}_j(f(.))\chi _{T_j}(.)\subset
{\rm Mod}\, {\mathcal F}(.;.)+{\rm Mod}\, f(.)\, .
$$
On the other hand,
$$
D^{(\rho _{\mathcal V}),W}\bigl( \, {\mathcal F}(f(.);.),\, \sum\limits_{j\in 
{\bf N}}\, {\mathcal F}_j(f(.))\chi _{T_j}(.)\, \bigr) <\epsilon \, .
$$
Hence, by Corollary \ref{c1}, we get ${\mathcal F}(f(.);.)\in W({\bf R},
{\mathcal V})$ and ${\rm Mod}\, {\mathcal F}(f(.);.)\subset {\rm Mod}\, {\mathcal 
F}(.;.)+{\rm Mod}\, f(.)$. 
\end{proof}

{\bf Remark 2}. From Lemma \ref{l9}, Theorems \ref{th1}, \ref{th2} and \ref{th4}
we obtain also the following assertion. Let $({\mathcal U},\rho )$ and 
$({\mathcal V},\rho _{\mathcal V})$ be complete metric spaces, let $r>0$ and let 
$p\geq 1$. Suppose that a function ${\bf R}\ni t\to {\mathcal F}(.;t)\in 
C({\mathcal U},{\mathcal V})$ satisfies the following two conditions:

(1) for any $x\in {\mathcal U}$ the function ${\bf R}\ni t\to {\mathcal F}
(.|_{B_r(x)};t)\in C(B_r(x),{\mathcal V})$ belongs to the space 
$$
W_1({\bf R},(C(B_r(x),{\mathcal V}),d_{C(B_r(x),{\mathcal V})}))
$$
(we denote by ${\mathcal F}(.|_Y)$ the restriction of a function ${\mathcal F}(.)
\in C({\mathcal U},{\mathcal V})$ to a non-empty set $Y\subset {\mathcal U}$);

(2) for a.e. $t\in {\bf R}$ the inequality 
$$
\rho _{\mathcal V}({\mathcal F}(x;t),y_0)\leq A\rho (x,x_0)+B(t)
$$
holds for all $x\in {\mathcal U}$, where $x_0\in {\mathcal U}$ and $y_0\in 
{\mathcal V}$ are some fixed points, $A\geq 0$, $B(.)\in M^{\sharp}_p({\bf R},
{\bf R})$.

Then for any function $f\in W_p({\bf R},{\mathcal U})$ we have ${\mathcal F}
(f(.);.)\in W_p({\bf R},{\mathcal V})$ and 
$$
{\rm Mod}\, {\mathcal F}(f(.);.)\subset {\rm Mod}\, f(.)+\sum\limits_{x\in 
{\mathcal U}}{\rm Mod}\, {\mathcal F}(.|_{B_r(x)};.)\, .
$$

\section{Proof of Theorem \ref{th7}}

\begin{lemma} \label{l11}
Let ${\mathbb F}\in {\mathcal A}^{(W)}$, $\Delta >0$. Then for any $\epsilon \in 
(0,1]$ there exist numbers $\delta =\delta (\epsilon ,\Delta )>0$, $l=l(\epsilon 
,\Delta ,{\mathbb F})>0$ and $\widetilde {\alpha}=\widetilde {\alpha}(\epsilon ,
\Delta ,{\mathbb F})>0$ such that for all $\alpha \geq \widetilde {\alpha}$ and 
all functions $f\in {\mathbb F}$
$$
\sup\limits_{\xi \in {\bf R}}\, {\rm meas}\, \{ t\in [\xi ,\xi +l]:|f(t)+\Delta
\sin \alpha t|<\delta \} < \epsilon l\, .
$$
\end{lemma}

\begin{proof}
Let us choose a number $N=N(\epsilon )\in {\bf N}$ for which $(N+1)^{-1}<\frac 
{\epsilon}3$ (then $N\geq 3$). Let $\epsilon ^{\, \prime}\doteq \frac 13\, 
\epsilon N^{-1}(N+1)^{-1}\leq \frac {\epsilon}{36}<1$, $\delta ^{\, \prime}=2\sin 
\frac {\pi}{2N}\, \sin \frac {\pi \epsilon ^{\, \prime}}6\, $, $\delta =\delta 
(\epsilon ,\Delta )=\min \, \{ 1,\frac 13\, \delta ^{\, \prime}\Delta \} $. There
are numbers $l=l(\epsilon ,\Delta ,{\mathbb F})>0$ and $\tau _0=\tau _0(\epsilon 
,\Delta ,{\mathbb F})\in (0,\frac 29\, \epsilon ]\subset (0,1)$ such that the
inequality 
$$
D^{(\rho )}_l(f(.),f(.+\tau ))<\epsilon ^{\, \prime}\delta 
$$
holds for all $f\in {\mathbb F}$ and $\tau \in [0,\tau _0]$. We define the number
$\widetilde {\alpha}=\pi \tau _0^{-1}$. Let $0<\tau \leq \tau_0\, $, $\alpha 
\doteq \pi \tau ^{-1}\geq \widetilde {\alpha}$. For $j=1,\dots ,N$ (and $f\in 
{\mathbb F}$, $\xi \in {\bf R}$) let us define the sets
$$
{\mathcal L}_j(\xi )={\mathcal L}_j(f,\tau ;\xi )=\{ t\in [\xi ,\xi +l]:|f(t+
\frac jN\, \tau )-f(t)|\geq \delta \} \, .
$$
We have
$$
{\rm meas}\, {\mathcal L}_j(\xi )\leq \frac 1{\delta}\, \int\limits_{\xi}^{\xi 
+l}\min \, \{ 1,|f(t+\frac jN\, \tau )-f(t)|\} \, dt<\epsilon ^{\, \prime}l\, .
$$
For $j=1,\dots ,N$ (and $\xi \in {\bf R}$) we also consider the sets  
$$
{\mathcal M}_j(\xi )={\mathcal M}_j(\tau ;\xi )=\{ t\in [\xi ,\xi +l]:|\cos
\alpha (t+\frac j{2N}\, \tau )\, \sin \frac {\alpha j}{2N}\, \tau \, |\leq \frac
{\delta ^{\, \prime}}2\, \} \, ,
$$ $$
{\mathcal M}_j^{\prime}(\xi )={\mathcal M}_j^{\prime}(\tau ;\xi )=\{ t\in [\xi ,
\xi +l]:|\cos (\alpha t+\frac {j\pi }{2N}\, )|\leq \frac 12\, \delta ^{\, \prime} 
\sin ^{-1}\frac {\pi}{2N}=\sin \frac {\pi \epsilon ^{\, \prime}}6\, \} \, ;
$$
${\mathcal M}_j(\xi )\subset {\mathcal M}_j^{\prime}(\xi )$ and ${\mathcal M}_j^
{\prime}(\xi )=[\xi ,\xi +l]\bigcap \, \bigl( \bigcup_{s\in {\bf Z}}\,
[\beta ^-_s\, ,\beta ^+_s]\bigr) $, where
$$
\beta ^{\pm}_s=\bigl( s+\frac 12\, \bigr) \, \frac {\pi}{\alpha}-\frac {j\pi }{2N
\alpha }\pm \frac {\pi \epsilon ^{\, \prime}}{6\alpha }\, .
$$
Let $\kappa $ be the number of closed intervals $[\beta ^-_s\, ,\beta ^+_s]$, 
$s\in {\bf Z}$, intersecting with the closed interval $[\xi ,\xi +l]$. We have 
the estimate $\kappa \leq (\frac {\pi}{\alpha})^{-1}(l+\frac {2\pi }{\alpha })
=\frac {l\alpha }{\pi}+2$, hence
$$
{\rm meas}\, {\mathcal M}_j(\xi )\leq {\rm meas}\, {\mathcal M}_j^{\prime}(\xi )
\leq \kappa \, \frac {\pi \epsilon ^{\, \prime}}{3\alpha }\leq \bigl( \, \frac l3
+\frac 23\, \tau _0\bigr) \, \epsilon ^{\, \prime}\leq l\epsilon ^{\, \prime}\, .
$$
In what follows, we suppose that the sets ${\mathcal L}_j(\xi )$ contain (in
addition) the numbers $t\in [\xi ,\xi +l]$ for which at the least one of the
functions $f(t)$, $f(t+\frac jN\, \tau )$ is not defined (these numbers form the
set of measure zero). Let
$$
{\mathcal L}(\xi )={\mathcal L}(f,\tau ;\xi )=\bigcup\limits_{j=1}^N\, \biggl(
\, \bigcup\limits_{s=0}^{N-j}({\mathcal L}_j(\xi )-\frac sN\, \tau )\biggr)
$$
(here ${\mathcal L}_j(\xi )-\frac sN\, \tau =\{ t=\eta -\frac sN\, \tau :\eta \in
{\mathcal L}_j(\xi )\} $). Since ${\rm meas}\, {\mathcal L}_j(\xi )<\epsilon ^{\, 
\prime}l$, $j=1,\dots ,N$, we get ${\rm meas}\, {\mathcal L}(\xi )<\frac 12
\, N(N+1)\epsilon ^{\, \prime}l=\frac 16\, \epsilon l$. If $t\in [\xi ,\xi +l-
\tau ]\backslash {\mathcal L}(\xi )$, then the numbers $t+\frac jN\, \tau $, $j=0,1,\dots ,N$, 
belong to the closed interval $[\xi ,\xi +l]$ and
$$
|f(t+\frac {j_1}N\, \tau )-f(t+\frac {j_2}N\, \tau )|<\delta 
$$
for all $j_1\, , j_2\in \{ 0,1,\dots ,N\} $. Let
$$
{\mathcal M}(\xi )={\mathcal M}(\tau ;\xi )=\bigcup\limits_{j=1}^N\, \biggl(
\, \bigcup\limits_{s=0}^{N-j}({\mathcal M}_j(\xi )-\frac sN\, \tau )\biggr) \, ;
$$
${\rm meas}\, {\mathcal M}(\xi )\leq \frac 12\, N(N+1)\epsilon ^{\, \prime}l=
\frac 16\, \epsilon l$. If $t\in [\xi ,\xi +l-\tau ]\backslash {\mathcal M}(\xi 
)$, then the numbers $t+\frac jN\, \tau $, $j=0,1,\dots ,N$, also belong to the 
closed interval $[\xi ,\xi +l]$ and for all $j_1\, , j_2\in \{ 0,1,\dots ,N\} $, 
$j_1<j_2\, $, we have
$$
|\Delta \sin \alpha (t+\frac {j_1}N\, \tau )-\Delta \sin \alpha (t+\frac {j_2}N\, 
\tau )|=
$$ $$
=2\Delta |\cos \alpha (t+\frac {j_1}N\, \tau +\frac {j_2-j_1}{2N}\, \tau )\, \sin
\alpha \, \frac {j_2-j_1}{2N}\, \tau |>\Delta \delta ^{\, \prime}\geq 3\delta \, .
$$
Let $G(t)=f(t)+\Delta \sin \alpha t$, $t\in {\bf R}$. We define (for $\xi \in
{\bf R}$) the set
$$
{\mathcal O}(\xi )={\mathcal O}(f,\tau ;\xi )=[\xi ,\xi +l-\tau ]\backslash
({\mathcal L}(\xi )\bigcup {\mathcal M}(\xi ))\, .
$$
For each $t\in {\mathcal O}(\xi )$ either $|G(t+\frac jN\, \tau )|\geq \delta $ 
for all $j=0,1,\dots ,N$ or there exists a number $j_0\in \{ 0,1,\dots ,N\} $ such
that $|G(t+\frac {j_0}N\, \tau )|< \delta $. Consider the minimal number $j_0$ 
for which the last inequality holds. If $j_0<N$, then for any $j\in \{ j_0+1,
\dots ,N\} $ we have
$$
|G(t+\frac jN\, \tau )-G(t+\frac {j_0}N\, \tau )|\geq
$$ $$
\geq |\Delta \sin \, \alpha (t+\frac jN\, \tau )- \Delta \sin \, \alpha (t+\frac 
{j_0}N\, \tau )|-|f(t+\frac jN\, \tau )-f(t+\frac {j_0}N\, \tau )|>3\delta -
\delta =2\delta \, ,
$$
and therefore, $|G(t+\frac jN\, \tau )|>\delta $. We have got that in the case 
$t\in {\mathcal O}(\xi )$ there is at most one number $t+\frac jN\, \tau $, $j=0,
1,\dots ,N$, such that $|G(t+\frac jN\, \tau )|<\delta $. Let denote by $\chi 
(t)$, $t\in {\bf R}$, the characteristic function of the set $\{ t\in [\xi ,\xi 
+l-\tau ]:|G(t)|<\delta \} $;
$$
\widetilde {\chi}(t)=\sum\limits_{j=0}^N\chi (t+\frac jN\, \tau )\, ,\ t\in {\bf
R}\, .
$$
Then $\widetilde {\chi}(t)\leq 1$ for all $t\in {\mathcal O}(\xi )$, and hence 
$$
\int\limits_{\xi}^{\xi +l-\tau }\widetilde {\chi}(t)\, dt=\int\limits_{{\mathcal
O}(\xi )}\widetilde {\chi}(t)\, dt+\int\limits_{[\, \xi ,\xi +l-\tau ]\, 
\backslash {\mathcal O}(\xi )}\widetilde {\chi}(t)\, dt\leq
$$ $$
\leq {\rm meas}\, {\mathcal O}(\xi )+(N+1)\, {\rm meas}\, ({\mathcal L}(\xi )
\bigcup {\mathcal M}(\xi ))<\bigl( 1+\frac 13\, (N+1)\epsilon \bigr) l\, .
$$
On the other hand,
$$
\int\limits_{\xi}^{\xi +l-\tau }\widetilde {\chi}(t)\, dt=(N+1)\int\limits_{\xi}
^{\xi +l-\tau }\chi (t)\, dt-\sum\limits_{j=1}^N\, \int\limits_{\xi}^{\xi +\frac
jN\, \tau }\chi (t)\, dt\geq
$$ $$
\geq (N+1)\, {\rm meas}\, \{ t\in [\xi ,\xi +l-\tau ]:|G(t)|<\delta \} -\frac 12\,
(N+1)\tau \geq
$$ $$
\geq (N+1)\, {\rm meas}\, \{ t\in [\xi ,\xi +l]:|G(t)|<\delta \} -\frac 32\,
(N+1)\tau \, .
$$
Hence
$$
{\rm meas}\, \{ t\in [\xi ,\xi +l]:|G(t)|<\delta \} <
$$ $$
<\biggl( \frac 1{N+1}+\frac {3\tau }{2l}+\frac {\epsilon}3\biggr) \, l\leq
\biggl( \epsilon -\biggl( \frac {\epsilon}3-\frac 1{N+1}\biggr) \biggr) \, l<
\epsilon l
$$
(for all $\xi \in {\bf R}$). 
\end{proof}

The following Lemma is an immediate consequence of Lemma \ref{l11}.

\begin{lemma} \label{l12}
Let ${\mathbb F}\in {\mathcal A}^{(W)}$, $\Delta >0$. Then for any $\epsilon \in 
(0,1]$ there exist numbers $\delta =\delta (\epsilon ,\Delta )>0$, $l=l(\epsilon 
,\Delta ,{\mathbb F})>0$ and $\widetilde {\alpha}=\widetilde {\alpha}(\epsilon ,
\Delta ,{\mathbb F})>0$ such that for each function $g\in L_{\infty}({\bf R},{\bf 
R})$ satisfying the condition $\| g\| _{L_{\infty}({\bf R},{\bf R})}\leq \delta 
$, and for all $\alpha \geq \widetilde {\alpha}$ and all functions $f\in {\mathbb 
F}$
$$
\sup\limits_{\xi \in {\bf R}}\, {\rm meas}\, \{ t\in [\xi ,\xi +l]:|f(t)+\Delta
\sin \alpha t+g(t)|<\delta \} < \epsilon l\, .
$$
\end{lemma}

{\it Proof of Theorem \ref{th7}}. Let $\Delta _0=\frac {\Delta}2\, $, $f_0(.)=
f(.)$ (for all functions $f\in {\mathbb F}$). By Lemma \ref{l12}, there are
numbers $\delta _0=\delta _0(\Delta )>0$, $l_0=l_0(\Delta ,{\mathbb F})>0$ and
$\alpha _0=\alpha _0(b,\Delta ,{\mathbb F})\in \frac {2\pi }b\, {\bf N}$ such that
for all functions $f_1(t)\doteq f_0(t)+\Delta _0\sin \alpha _0t$, $t\in {\bf R}$, 
and all functions $\widetilde g_1\in L_{\infty}({\bf R},{\bf R})$ that satisfy
the condition $\| \widetilde g_1\| _{L_{\infty}({\bf R},{\bf R})}\leq \delta _0\, 
$, the inequality
$$
\sup\limits_{\xi \in {\bf R}}\, {\rm meas}\, \{ t\in [\xi ,\xi +l_0]:|f_1(t)+
\widetilde g_1(t)|<\delta _0\} <2^{-1}l_0 \eqno (3)
$$
holds, furthermore $\{ f_1(.):f\in {\mathbb F}\} \in {\mathcal A}^{(W)}$. We 
shall successively for $j=1,2,\dots $ find numbers $\Delta _j=\Delta _j(\Delta )
>0$, $\delta _j=\delta _j(\Delta )>0$, $l_j=l_j(\Delta ,{\mathbb F})>0$,
$\alpha _j=\alpha _j(b,\Delta ,{\mathbb F})\in \frac {2\pi }b\, {\bf N}$ and
functions $f_{j+1}\in W({\bf R},{\bf R})$ dependent on $f_j\, $, $\Delta _j\, $ 
and $\alpha _j\, $, for which $\{ f_{j+1}(.):f\in {\mathbb F}\} \in {\mathcal A}
^{(W)}$. If the numbers $\Delta _k\, $, $\delta _k\, $, $l_k\, $, $\alpha _k$ and
the functions $f_{k+1}$ have been found for all $k=0,\dots ,\, j-1$, where $j\in 
{\bf N}$, then we choose the number $\Delta _j=\Delta _j(\Delta )>0$ such that 
the inequalities $\Delta _j<2^{-(j+1)}\Delta $, $\Delta _j\leq 2^{-j}\delta _0\, 
$, $\Delta _j\leq 2^{-(j-1)}\delta _1\, $, $\dots $ , $\Delta _j\leq 2^{-1}\delta 
_{j-1}\, $ hold. Further (according to Lemma \ref{l12}), choose numbers $\delta 
_j=\delta _j(\Delta )>0$, $l_j=l_j(\Delta ,{\mathbb F})>0$ and $\alpha _j=\alpha 
_j(b,\Delta ,{\mathbb F})\in \frac {2\pi }b\, {\bf N}$ such that for all functions 
$f_{j+1}(t)\doteq f_j(t)+\Delta _j\sin \alpha _jt$, $t\in {\bf R}$, and all
functions $\widetilde g_{j+1}\in L_{\infty}({\bf R},{\bf R})$ satisfying the
condition $\| \widetilde g_{j+1}\| _{L_{\infty}({\bf R},{\bf R})}\leq \delta _j\, 
$, the inequality
$$
\sup\limits_{\xi \in {\bf R}}\, {\rm meas}\, \{ t\in [\xi ,\xi +l_j]:|f_{j+1}(t)+
\widetilde g_{j+1}(t)|<\delta _j\} <2^{-j-1}l_j \eqno (4)
$$
holds. We also have $\{ f_{j+1}(.):f\in {\mathbb F}\} \in {\mathcal A}^{(W)}$. 
Next, let us set
$$
g(t)=\sum\limits_{j=0}^{+\infty}\Delta _j\sin \, \alpha _jt\, ,\ t\in {\bf R}\, .
$$
Since $\Delta _0=\frac {\Delta}2$ and $\Delta _j<2^{-(j+1)}\Delta $ for all $j\in 
{\bf N}$, it follows that the function $g(.)$ is continuous and $b$-periodic,
moreover,
$$
\| g\| _{L_{\infty}({\bf R},{\bf R})}\leq \sum\limits_{j=0}^{+\infty}\Delta _j
<\Delta \, .
$$
We define the functions
$$
g_j(t)=\sum\limits_{k=j}^{+\infty}\Delta _k\sin \, \alpha _kt\, ,\ t\in {\bf R}\, 
,\, j\in {\bf N}\, .
$$ 
For all $t\in {\bf R}$ we have
$$
|g_j(t)|\leq \sum\limits_{k=j}^{+\infty}\Delta _k\leq \bigl( \frac 12+\frac 14+
\dots \bigr) \, \delta _{j-1}=\delta _{j-1}\, ,
$$
hence, it follows from (3) and (4) that for all numbers $j=0,1,\dots $ (and all
functions $f\in {\mathbb F}$) the inequality
$$
\sup\limits_{\xi \in {\bf R}}\, {\rm meas}\, \{ t\in [\xi ,\xi +l_j]:|f(t)+g(t)|
<\delta _j\} <2^{-j-1}l_j\, .  
$$
holds. The proof of Theorem \ref{th7} is complete.

\section{Proof of Theorem \ref{th3}}

\begin{theorem} \label{th8} 
Let $({\mathcal U},\rho )$ be a complete metric space, let $F\in W({\bf R},{\rm 
cl}\, {\mathcal U})$ and let $g\in W({\bf R},{\mathcal U})$. Then for any 
$\epsilon >0$ there exists a function $f\in W({\bf R},{\mathcal U})$ such that 
${\rm Mod}\, f\subset {\rm Mod}\, F+{\rm Mod}\, g\, $, $f(t)\in F(t)$ a.e. and 
$\rho (f(t),g(t))<\rho (g(t),F(t))+\epsilon $ a.e.
\end{theorem}

\begin{proof} 
Let number $\epsilon \in (0,1]$ be fixed. We choose numbers $\gamma _n>0$, $n\in 
{\bf N}$, such that
$$
\sum\limits_{n=1}^{+\infty}(\gamma _n+\gamma _{n+1})<\frac 16
$$
(then $\gamma _1<\frac 16\, $). From Theorem \ref{th2}, Lemma \ref{l5} and
Corollary \ref{c3} it follows that for each $n\in {\bf N}$ there exist sets 
$F^{(n)}_j\in {\rm cl}\, {\mathcal U}$, points $g^n_j\in {\mathcal U}$ and 
disjoint measurable (in the Lebesgue sense) sets $T^{(n)}_j\subset {\bf R}$, 
$j\in {\bf N}$, such that $\{ T^{(n)}_j\} _{j\in {\bf N}}\in {\mathfrak M}^{(W)}
({\rm Mod}\, F+{\rm Mod}\, g)$, the functions $F(t)$ and $g(t)$ are defined for
all $t\in \bigcup_jT^{(n)}_j$, and for all $t\in T^{(n)}_j$, $j\in {\bf N}$, 
we have ${\rm dist}_{\rho ^{\, \prime}}(F(t),F^{(n)}_j)<\gamma _n\epsilon <1$ and 
$\rho (g(t),g^n_j)<\gamma _n\epsilon \, $. Let 
$$
T=\bigcap\limits_n\bigcup\limits_jT^{(n)}_j\, ;
$$
${\rm meas}\, {\bf R}\backslash T=0$. By Corollary \ref{c3}, for every $n\in 
{\bf N}$
$$
\{ T^{(1)}_{j_1}\bigcap \dots \bigcap T^{(n)}_{j_n}\} _{j_s\in {\bf N},\ s=1,
\dots ,n}\in {\mathfrak M}^{(W)}({\rm Mod}\, F+{\rm Mod}\, g)\, .
$$
With each number $n\in {\bf N}$ and each collection $\{ j_1\, ,\dots ,\, j_n\} $ 
of indices $j_s\in {\bf N}$, $s=1,\dots ,n$, if $T^{(1)}_{j_1}\bigcap \dots 
\bigcap T^{(n)}_{j_n}\neq \emptyset $, we associate some point $f_{j_1\dots j_n}
\in F^{(n)}_{j_n}\subset {\mathcal U}$. These points are determined successively
for $n=1,2,\dots $. For $n=1$ we choose points $f_{j_1}\in F^{(1)}_{j_1}$ such
that the inequalities
$$
\rho (f_{j_1}\, ,g^1_{j_1})<\frac {\epsilon}6+\rho (g^1_{j_1}\, ,F^{(1)}_{j_1})
$$
hold. If points $f_{j_1\dots j_{n-1}}\in F^{(n-1)}_{j_{n-1}}$ have been found for
some $n\geq 2$, then we choose points $f_{j_1\dots j_{n-1}j_n}\in F^{(n)}_{j_n}$ 
such that
$$
\rho (f_{j_1\dots j_{n-1}}\, ,f_{j_1\dots j_{n-1}j_n})=\rho ^{\, \prime}(f_{j_1
\dots j_{n-1}}\, ,f_{j_1\dots j_{n-1}j_n})\leq  \eqno (5)
$$ $$
\leq 2\, {\rm dist}_{\rho ^{\, \prime}}(F^{(n-1)}_{j_{n-1}},F^{(n)}_{j_n})<2\,
(\gamma _{n-1}+\gamma _n)\epsilon <\frac {\epsilon}3\leq \frac 13\, .
$$
Let us define functions
$$
f(n;t)=\sum\limits_{j_1\, ,\dots ,\, j_n}f_{j_1\dots j_n}\chi _{T^{(1)}_{j_1}
\bigcap \dots \bigcap T^{(n)}_{j_n}}(t)\, ,\ t\in T\, ,\ n\in {\bf N}\, .
$$
According to Theorem \ref{th1}, we have $f(n;.)\in W({\bf R},{\mathcal U})$ and 
${\rm Mod}\, f(n;.)\subset {\rm Mod}\, F+{\rm Mod}\, g\, $. It follows from (5) 
that the inequality
$$
\rho (f(n-1;t),f(n;t))< 2\, (\gamma _{n-1}+\gamma _n)\epsilon \eqno (6)
$$
holds for all $t\in T$ and $n\geq 2$. Since the metric space ${\mathcal U}$ is
complete, we obtain from (6) that the sequence of functions $f(n;.)$, $n\in {\bf 
N}$, converges as $n\to +\infty $ uniformly on the set $T\subset {\bf R}$ 
(therefore, in the metric $D^{(\rho ),W}$ as well) to a function $f(.)\in W({\bf 
R},{\mathcal U})$ for which ${\rm Mod}\, f\subset \sum_n{\rm Mod}\, f(n;.)\subset 
{\rm Mod}\, F+{\rm Mod}\, g\, $. We have $f(n;t)\in F^{(n)}_{j_n}$ and ${\rm 
dist}_{\rho ^{\, \prime}}(F(t),F^{(n)}_{j_n})<\gamma _n\epsilon <\frac 16$ for 
all $t\in T^{(n)}_{j_n}\bigcap T$. Since $\gamma _n\to 0$ as $n\to +\infty $, it
follows from this that $f(t)\in F(t)$ for all $t\in T$ (for a.e. $t\in {\bf R}$). 
With each number $t\in T$ we associate an infinite collection of indices $\{ j_1
\, ,\dots ,\, j_n\, ,\dots \} $ in such a way that $t\in T^{(n)}_{j_n}$, $n\in 
{\bf N}$. Then (for all $t\in T$)
$$
\rho (f(t),g(t))\leq \sum\limits_{n=1}^{+\infty}\rho (f_{j_1\dots j_n}\, ,f_{j_1
\dots j_nj_{n+1}})+\rho (f_{j_1}\, ,g^1_{j_1})+\rho (g^1_{j_1}\, ,g(t))<
$$ $$
<2\, \sum\limits_{n=1}^{+\infty}(\gamma _n+\gamma _{n+1})\epsilon +\frac 
{\epsilon}3+\rho (g^1_{j_1}\, ,F^{(1)}_{j_1})<
$$ $$
<\frac {2\epsilon }3+|\rho (g^1_{j_1}\, ,F^{(1)}_{j_1})-\rho (g^1_{j_1}\, 
,F(t))|+|\rho (g^1_{j_1}\, ,F(t))-\rho (g(t),F(t))|+\rho (g(t),F(t))<
$$ $$
<\frac {2\epsilon }3+\gamma _1\epsilon +\gamma _1\epsilon +\rho (g(t),F(t))<
\epsilon +\rho (g(t),F(t))\, .
$$
\end{proof}

{\bf Remark 3}. If conditions of Theorem \ref{th8} are fulfilled and, moreover,
$F\in W_p({\bf R},{\rm cl}_b\, {\mathcal U})\subset W({\bf R},{\rm cl}\, 
{\mathcal U})$, $p\geq 1$, then it follows from Theorem \ref{th4} that $f\in 
W_p({\bf R},{\mathcal U})$. Indeed, for a.e. $t\in {\bf R}$ we have
$$
\rho (x_0\, ,f(t))\leq \sup\limits_{x\in F(t)}\rho (x_0\, ,x)={\rm dist}\, (\{ 
x_0\} ,F(t))\, ,
$$
furthermore, ${\rm dist}\, (\{ x_0\} ,F(.))\in M^{\sharp}_p({\bf R},{\bf R})$. 
Hence $f(.)\in M^{\sharp}_p({\bf R},{\mathcal U})\bigcap W({\bf R},{\mathcal U})
=W_p({\bf R},{\mathcal U})$.

\begin{corollary} \label{c6}
Let $({\mathcal U},\rho )$ be a complete separable metric space and let $F\in W({\bf R},
{\rm cl}\, {\mathcal U})$. Then there exist functions $f_j\in W({\bf R},{\mathcal 
U})$, $j\in {\bf N}$, such that ${\rm Mod}\, f_j\subset {\rm Mod}\, F$ and $F(t)
=\overline {\bigcup_jf_j(t)}$ for a.e. $t\in {\bf R}$ (if $F\in W_p({\bf R},{\rm 
cl}_b\, {\mathcal U})\subset W({\bf R},{\rm cl}\, {\mathcal U})$, $p\geq 1$, then
all functions $f_j$ belong to the space $W_p({\bf R},{\mathcal U})$ (see Remark
3)).
\end{corollary}

\begin{proof} 
Let us choose points $x_k\in {\mathcal U}$, $k\in {\bf N}$, which form a 
countable dense set of the metric space ${\mathcal U}$. By Theorem \ref{th8},
for all $k,n\in {\bf N}$ there are functions $f_{k,n}\in W({\bf R},{\mathcal 
U})$ such that ${\rm Mod}\, f_{k,n}\subset {\rm Mod}\, F$, $f_{k,n}(t)\in F(t)$ 
a.e. and $\rho (f_{k,n}(t),x_k)<2^{-n}+\rho (x_k\, ,F(t))$ a.e. It remains to
renumber the functions $f_{k,n}(.)$ by a single index $j\in {\bf N}$.
\end{proof}

{\it Proof of Theorem \ref{th3}}. It can be assumed without loss of generality
that $\eta (t)\to 0$ as $t\to +0$. Since $F\in W({\bf R},{\rm cl}\, {\mathcal 
U})$ and $g\in W({\bf R},{\mathcal U})$, it follows that $\rho (g(.),F(.))\in 
W({\bf R},{\bf R})$ and ${\rm Mod}\, \rho (g(.),F(.))\subset {\rm Mod}\, F+{\rm 
Mod}\, g\, $. According to Corollary \ref{c5}, for each $j\in {\bf N}$ we choose
sets $T_j\in W({\bf R})$ such that ${\rm Mod}\, T_j\subset {\rm Mod}\, \rho 
(g(.),F(.))\subset {\rm Mod}\, F+{\rm Mod}\, g\, $, $\rho (g(t),F(t))<2^{-j}$ for
all $t\in T_j\, $, and $\rho (g(t),F(t))>2^{-j-1}$ for a.e. $t\in {\bf R}
\backslash T_j\, $. Further (after deletion of some subsets of measure zero from
the sets $T_j\, $, $j\in {\bf N}$), we can assume that $T_{j+1}\subset T_j\, $. 
Let $T_0={\bf R}$. We have $T_{j-1}\backslash T_j\in W({\bf R})$ and ${\rm Mod}\, 
T_{j-1}\backslash T_j\subset {\rm Mod}\, F+{\rm Mod}\, g\, $, $j\in {\bf N}$ (see
Lemma \ref{l5}). For each $j\in {\bf N}$, according to Theorem \ref{th8}, we
choose functions $f_j\in W({\bf R},{\mathcal U})$ for which ${\rm Mod}\, f_j
\subset {\rm Mod}\, F+{\rm Mod}\, g\, $, $f_j(t)\in F(t)$ a.e. and $\rho (f_j(t),
g(t))<\rho (g(t),F(t))+\eta (2^{-j-1})$ a.e. We define the functions
$$
f(.)=\sum\limits_{j=1}^{+\infty}f_j(.)\chi _{T_{j-1}\backslash T_j}(.)+g(.){\chi}
_{\bigcap\limits_jT_j}(.)\, ,
$$ $$
f(n;.)=\sum\limits_{j=1}^nf_j(.)\chi _{T_{j-1}\backslash T_j}(.)+f_{n+1}(.){\chi}
_{T_n}(.)\, ,\ n\in {\bf N}\, .
$$
By Theorem \ref{th1}, $f(n;.)\in W({\bf R},{\mathcal U})$ and ${\rm Mod}\, f(n;.)
\subset {\rm Mod}\, F+{\rm Mod}\, g\, $. Since $f_j(t)\in F(t)$ a.e. and $g(t)\in 
F(t)$ for all $t\in \bigcap_jT_j\, $, it follows that $f(n;t)\in F(t)$ and $f(t)
\in F(t)$ a.e. as well. For each $n\in {\bf N}$ we have 
$$
\rho (f(t),f(n;t))=0
$$ 
for a.e. $t\in {\bf R}\backslash T_{n+1}\, $, 
$$
\rho (f(t),f(n;t))=\rho (f_{m+1}(t),f_{n+1}(t))\leq 
\rho (f_{m+1}(t),g(t))+\rho (f_{n+1}(t),g(t))<
$$ $$
<2\, \rho (g(t),F(t))+\eta (2^{-m-2})+\eta (2^{-n-2})<2^{-n}+2\, \eta (2^{-n-2})
$$ 
for a.e. $t\in T_m\backslash T_{m+1}\, $, $m\geq n+1$, and 
$$
\rho (f(t),f(n;t))=\rho (g(t),f_{n+1}(t))<\eta (2^{-n-2})
$$ 
for a.e. $t\in \bigcap_jT_j\, $. Therefore
$$
{\rm ess}\, \sup\limits_{\hskip -0.7cm t\in {\bf R}}\, \rho (f(t),f(n;t))\to 0
$$
as $n\to +\infty $, hence $f\in W({\bf R},{\mathcal U})$ and ${\rm Mod}\, f
\subset \sum_n {\rm Mod}\, f(n;.)\subset {\rm Mod}\, F+{\rm Mod}\, g\, $. For 
a.e. $t\in T_{j-1}\backslash T_j\, $, $j\in {\bf N}$, the estimate
$$
\rho (f(t),g(t))=\rho (f_j(t),g(t))<\rho (g(t),F(t))+\eta (2^{-j-1})\leq
$$ $$
\leq \rho (g(t),F(t))+\eta (\rho (g(t),F(t)))
$$
holds. From this (since $f(t)=g(t)\in F(t)$ for $t\in \bigcap_jT_j$) we obtain
that for a.e. $t\in {\bf R}$
$$
\rho (f(t),g(t))\leq \rho (g(t),F(t))+\eta (\rho (g(t),F(t)))\, .
$$
If $F\in W_p({\bf R},{\rm cl}_b\, {\mathcal U})$, $p\geq 1$, then also $f\in 
W_p({\bf R},{\mathcal U})$ (see Remark 3). \hfill $\square$

The following Theorems can be proved (using Theorems \ref{th1} and \ref{th2},
Lemma \ref{l5} and Corollaries \ref{c2}, \ref{c3} and \ref{c5}) by analogy
with appropriate assertions on Stepanov a.p. functions and multivalued maps
\cite{D6,D7}.

For non-empty set $F\subset {\mathcal U}$ we use the notation
$F^{\, \epsilon}=\{ x\in {\mathcal U}:\rho (x,F)<\epsilon \} $, $\epsilon >0$.

The points $x_j\in {\mathcal U}$, $j=1,\dots ,n$, are said to form {\it
$\epsilon$-net} for (non-empty) set $F\subset {\mathcal U}$, $\epsilon >0$, if
$F\subset \bigl( \bigcup_jx_j\bigr) ^{\epsilon}$.

\begin{theorem} \label{th9}
Let $({\mathcal U},\rho )$ be a complete metric space, let $F\in W({\bf R},{\rm 
cl}_b\, {\mathcal U})$ and let $\epsilon >0$, $n\in {\bf N}$. Suppose that for
a.e. $t\in {\bf R}$ there are points $x_j(t)\in F(t)$, $j=1,\dots ,n$, which
form $\epsilon$-net for the set $F(t)$. Then for any $\epsilon ^{\, \prime}>
\epsilon $ there exist functions $f_j\in W({\bf R},{\mathcal U})$, $j=1,\dots 
,n$, such that ${\rm Mod}\, f_j\subset {\rm Mod}\, F$, $f_j(t)\in F(t)$ a.e.
and for a.e. $t\in {\bf R}$ the points $f_j(t)$, $j=1,\dots ,n$, form $\epsilon
^{\, \prime}$-net for the set $F(t)$.
\end{theorem}

\begin{corollary} \label{c7}
Let $({\mathcal U},\rho )$ be a compact metric space. Then a multivalued map
${\bf R}\ni t\to F(t)\in {\rm cl}\, {\mathcal U}={\rm cl}_b\, {\mathcal U}$
belongs to the space $W({\bf R},{\rm cl}\, {\mathcal U})=W_1({\bf R},{\rm 
cl}_b\, {\mathcal U})$ if and only if for each $\epsilon >0$ there exist a
number $n\in {\bf N}$ and functions $f_j\in W({\bf R},{\mathcal U})=W_1({\bf R},
{\mathcal U})$, $j=1,\dots ,n$, such that $f_j(t)\in F(t)$ a.e. and points 
$f_j(t)$, $j=1,\dots ,n$, for a.e. $t\in {\bf R}$ form $\epsilon$-net for the set 
$F(t)$ (furthermore, the functions $f_j$ for the multivalued map $F\in W({\bf R},
{\rm cl}\, {\mathcal U})$ can be chosen in such a way that ${\rm Mod}\, f_j
\subset {\rm Mod}\, F$).
\end{corollary}

\begin{theorem} \label{th10}
Let $({\mathcal U},\rho )$ be a compact metric space. Then a multivalued map
${\bf R}\ni t\to F(t)\in {\rm cl}\, {\mathcal U}$ belongs to the space $W({\bf R},
{\rm cl}\, {\mathcal U})$ if and only if there exist functions $f_j\in W({\bf R},
{\mathcal U})$, $j\in {\bf N}$, such that $F(t)=\overline {\bigcup_jf(t)}$ 
a.e. and the set $\{ f_j(.): j\in {\bf N}\} $ is precompact in the metric space
$L_{\infty}({\bf R},{\mathcal U})$ (furthermore, the functions $f_j$ for the 
multivalued map $F\in W({\bf R},{\rm cl}\, {\mathcal U})$ can be chosen in such 
a way that ${\rm Mod}\, f_j\subset {\rm Mod}\, F$).
\end{theorem}

\begin{theorem} \label{th11}
Let $({\mathcal U},\rho )$ be a complete metric space, let $F\in W({\bf R},{\rm 
cl}_b\, {\mathcal U})$, $\epsilon >0$, $\delta >0$, $n\in {\bf N}$, and let 
$g_j\in W({\bf R},{\mathcal U})$, $j=1,\dots ,n$. Suppose that for
a.e. $t\in {\bf R}$ the set of points $x_j(t)=g_j(t)$, for which $g_j(t)\in
(F(t))^{\delta}$, can be supplemented (if it consists of less than $n$ points)
to $n$ points $x_j(t)\in (F(t))^{\delta}$, $j=1,\dots ,n$, which
form $\epsilon$-net for the set $F(t)$ (coincident points with different
indices are considered here as different points). Then for any $\epsilon ^{\, \prime}>
\epsilon +\delta $ there exist functions $f_j\in W({\bf R},{\mathcal U})$, $j=1,
\dots ,n$, such that ${\rm Mod}\, f_j\subset {\rm Mod}\, F+\sum_{k=1}^n
{\rm Mod}\, g_k\, $, $f_j(t)\in F(t)$ a.e., $f_j(t)=g_j(t)$ for a.e. $t\in
\{ \tau \in {\bf R}:g_j(\tau )\in F(\tau )\} $ and the points $f_j(t)$, $j=1,
\dots ,n$, for a.e. $t\in {\bf R}$ form $\epsilon^{\, \prime}$-net for the set $F(t)$.
\end{theorem}

\end{document}